\documentclass{article}

\usepackage{amsmath}
\usepackage{amssymb}

\newtheorem{thm}{Theorem}
\newtheorem{conj}{Conjecture}
\newtheorem{lem}{Lemma}
\newtheorem{prop}{Proposition}

\newtheorem{ques}{Question}
\newtheorem{defn}{Definition}

\title{Finding Almost Squares III}
\author{Tsz Ho Chan}

\begin{document}
\maketitle
\begin{abstract}
An almost square of type $2$ is an integer $n$ that can be factored
in two different ways as $n = a_1 b_1 = a_2 b_2$ with $a_1$, $a_2$,
$b_1$, $b_2 \approx \sqrt{n}$. In this paper, we shall improve upon
previous result on short intervals containing an almost square of
type $2$. This leads to an inquiry of finding a short interval
around $x$ that contains an integer divisible by some integer in
$[x^c, 2 x^c]$ with $0 < c < 1$.
\end{abstract}
\section{Introduction and main results}
In [\ref{C1}] and [\ref{C2}], the author started studying almost
square, an integer $n$ that can be factored as $n = a b$ with $a ,
b$ close to $\sqrt{n}$. More specifically, for $0 \leq \theta \leq
1/2$ and $C > 0$,
\begin{defn}
An integer $n$ is a ($\theta$, $C$)-almost square of type $1$ if $n
= a b$ for some integers $a, b$ in the interval $[n^{1/2} - C
n^\theta, n^{1/2} + C n^\theta]$.
\end{defn}
\begin{defn}
An integer $n$ is a ($\theta$, $C$)-almost square of type $2$ if $n
= a_1 b_1 = a_2 b_2$ for some integers $a_1 < a_2 \leq b_2 < b_1$ in
the interval $[n^{1/2} - C n^\theta, n^{1/2} + C n^\theta]$.
\end{defn}
Let $x$ be a large positive real number. Following [\ref{C1}] and
[\ref{C2}], we are interested in finding almost squares of type $1$
or $2$ near to $x$. In particular, given $0 \leq \theta \leq 1/2$,
we want to find ``admissible" $\phi_i \geq 0$ (as small as possible)
such that, for some constants $C_{\theta, i}, D_{\theta, i} > 0$,
the interval $[x - D_{\theta, i} x^{\phi_i}, x + D_{\theta, i}
x^{\phi_i}]$ contains a ($\theta$, $C_{\theta, i}$)-almost square of
type $i$ ($i = 1,2$) for all large $x$.
\begin{defn}
$$f(\theta) := \inf \phi_1 \; \hbox{ and } \; g(\theta) := \inf
\phi_2$$ where the infima are taken over all the ``admissible"
$\phi_i$ ($i = 1, 2$) respectively.
\end{defn}
Clearly $f$ and $g$ are non-increasing functions of $\theta$.
Summarizing the results in [\ref{C1}] and [\ref{C2}], we have
\begin{thm} \label{thm1}
For $0 \leq \theta \leq 1/2$,
\[
f(\theta) \left\{
\begin{tabular}{ll}
$= 1/2$, & if $0 \leq \theta < 1/4$, \\
$= 1/4$, & if $\theta = 1/4$, \\
$= 1/2 - \theta$, & if $1/4 \leq \theta \leq 3/10$ and a conjectural
upper bound on \\
 & certain average of twisted incomplete Salie sum is true, \\
$\geq 1/2 - \theta$, & if $1/4 \leq \theta \leq 1/2$.
\end{tabular} \right.
\]
\end{thm}
\begin{thm} \label{thm2}
For $0 \leq \theta \leq 1/2$,
\[
g(\theta) \left\{
\begin{tabular}{ll}
does not exist, & if $0 \leq \theta < 1/4$, \\
$\geq 1 - 2\theta$, & if $1/4 \leq \theta \leq 1/2$, \\
$\leq 1 - \theta$, & if $1/4 \leq \theta \leq 1/3$.
\end{tabular} \right.
\]
\end{thm}
And we conjectured that
\begin{conj}
For $0 \leq \theta \leq 1/2$,
\[
f(\theta) = \left\{
\begin{tabular}{ll}
$1/2$, & if $0 \leq \theta < 1/4$, \\
$1/2 - \theta$, & if $1/4 \leq \theta \leq 1/2$;
\end{tabular} \right.
\]
and
\[
g(\theta) = \left\{
\begin{tabular}{ll}
does not exist, & if $0 \leq \theta < 1/4$, \\
$1 - 2\theta$, & if $1/4 \leq \theta \leq 1/2$.
\end{tabular} \right.
\]
\end{conj}

In this paper, we are going to improve Theorem \ref{thm2}.
\begin{thm} \label{thm3}
For $1/4 \leq \theta \leq 1/2$,

\begin{tabular}{lll}
(i) & $g(1/4) \leq 5/8$, & \\
(ii) & $g(\theta) \leq 9/16$, & if $5/16 \leq \theta \leq 1/2$, \\
(iii) & $g(\theta) \leq 17/32$, & if $5/16 \leq \theta \leq 1/2$, \\
(iv) & $g(\theta) \leq 1/2$, & if $1/3 < \theta \leq 1/2$, \\
(v) & $g(\theta) \leq 1/2$, & if $743/2306 < \theta \leq 1/2$.
\end{tabular}
\end{thm}
Clearly (iii) is better than (ii). The reason we keep (ii) is that
(ii) and (iii) use different approaches. Also (v) includes (iv). The
reason we keep (iv) is that (iv) provides a prototype for (v).

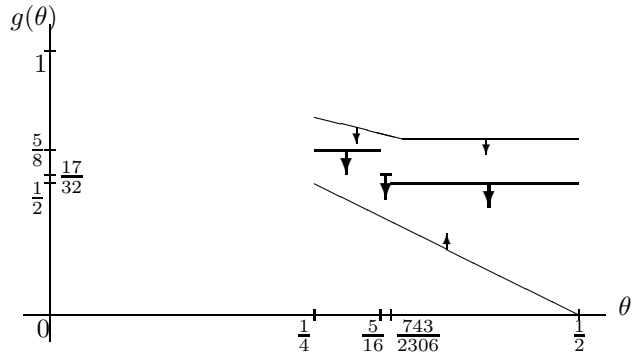
\begin{picture}(220,140)(-10,-10)
\put(-10,0){\line(1,0){220}} \put(0,-10){\line(0,1){120}}
\put(-2,50){\line(1,0){4}} \put(-2,100){\line(1,0){4}}
\put(100,-2){\line(0,1){4}} \put(200,-2){\line(0,1){4}}
\put(-5,-8){$0$} \put(-8,42){$\frac{1}{2}$} \put(-6,92){$1$}
\put(-15,110){$g(\theta)$} \put(93,-10){$\frac{1}{4}$}
\put(197,-10){$\frac{1}{2}$} \put(215,0){$\theta$}
\put(125,-2){\line(0,1){4}} \put(117,-10){$\frac{5}{16}$}
\put(129,-2){\line(0,1){4}} \put(130,-10){$\frac{743}{2306}$}
\put(-2,62.5){\line(1,0){4}} \put(-8,60){$\frac{5}{8}$}
\put(-2,53.125){\line(1,0){4}} \put(3,50){$\frac{17}{32}$}
\thicklines \put(100,62.5){\line(1,0){25}}
\put(125,53.125){\line(1,0){4}} \put(129,50){\line(1,0){71}}
\put(112,62.5){\vector(0,-1){9}} \put(127,53){\vector(0,-1){9}}
\put(166,50){\vector(0,-1){9}} \thinlines
\put(100,50){\line(2,-1){100}} \put(150,25){\vector(0,1){6}}
\put(100,75){\line(4,-1){33.3}} \put(116,71){\vector(0,-1){6}}
\put(133.3,66.675){\line(1,0){66.7}}
\put(165,66.675){\vector(0,-1){6}}
\end{picture}

\bigskip
The above picture summarizes Theorems \ref{thm2} and \ref{thm3}. The
thin line segments are the upper and lower bounds from Theorem
\ref{thm2}. The thick line segments are the upper bounds from
Theorem \ref{thm3}. The next challenge is to beat the $\frac{1}{2}$
upper bound for $g(\theta)$.

\bigskip

{\bf Some Notations:} Throughout the paper, $\epsilon$ denotes a
small positive number. Both $f(x) = O(g(x))$ and $f(x) \ll g(x)$
mean that $|f(x)| \leq C g(x)$ for some constant $C > 0$. Moreover
$f(x) = O_\lambda(g(x))$ and $f(x) \ll_\lambda g(x)$ mean that the
implicit constant $C = C_\lambda$ may depend on the parameter
$\lambda$. Finally $f(x) \asymp g(x)$ means that $f(x) \ll g(x)$ and
$g(x) \ll f(x)$.
\section{Proof of Theorem \ref{thm3} (i)}

Let $1/4 \leq \theta \leq 1/2$. From [\ref{C2}], we recall that a
($\theta$, $C$)-almost square of type $2$ must be of the form
$$n = (d_1 e_1) (d_2 e_2) = (d_1 e_2) (d_2 e_1)$$
where $a_1 = d_1 e_1$, $b_1 = d_2 e_2$, $a_2 = d_1 e_2$, $b_2 = d_2
e_1$; $n^{1/2} - C n^\theta \leq a_1 < a_2 \leq b_2 < b_1 \leq
n^{1/2} + C n^\theta$;
$$\frac{1}{2 C} n^{\frac{1}{2} - \theta} - \frac{1}{2} \leq d_1,
d_2, e_1, e_2 \leq 2C n^\theta, \; e_2 - e_1 \leq 2C
\frac{n^\theta}{d_2}, \; d_2 - d_1 \leq 2C \frac{n^\theta}{e_2}.$$
Let $1 \leq k \ll 1$ be any integer. By $\theta = 1/4$ case in
Theorem \ref{thm1}, for some constant $C > 0$, we can find integers
$d, e \in [x^{1/4} - C x^{1/8}, x^{1/4} + C x^{1/8}]$ such that
$$d e = x^{1/2} - 2k x^{1/4} + O(x^{1/8}).$$
Then
$$(d + 2k) (e + 2k) = d e + 2k (d+e) + k^2 = x^{1/2} + 2k x^{1/4}
+ O(x^{1/8}),$$ and
$$d e (d + 2k) (e + 2k) = x - 4k^2 x^{1/2} + O(x^{5/8}) = x +
O(x^{5/8}).$$ This gives $g(1/4) \leq 5/8$.
\section{Proof of Theorem \ref{thm3} (ii)}
The key idea is the identity
$$a b = \Bigl(\frac{a+b}{2}\Bigr)^2 - \Bigl(\frac{a-b}{2}\Bigr)^2$$
as used in [\ref{C1}]. Using this identity,
\begin{align*}
d_1 e_1 d_2 e_2 =& \Bigl[\Bigl(\frac{d_2 + d_1}{2}\Bigr)^2 -
\Bigl(\frac{d_2 - d_1}{2}\Bigr)^2 \Bigr] \Bigl[\Bigl(\frac{e_2 +
e_1}{2}\Bigr)^2 - \Bigl(\frac{e_2 - e_1}{2}\Bigr)^2 \Bigr] \\
=& \Bigl(\frac{d_2 + d_1}{2}\Bigr)^2 \Bigl(\frac{e_2 +
e_1}{2}\Bigr)^2 - \Bigl(\frac{d_2 - d_1}{2}\Bigr)^2 \Bigl(\frac{e_2
+ e_1}{2}\Bigr)^2 \\
&- \Bigl(\frac{e_2 - e_1}{2}\Bigr)^2 \Bigl(\frac{d_2 +
d_1}{2}\Bigr)^2 + \Bigl(\frac{d_2 - d_1}{2}\Bigr)^2 \Bigl(\frac{e_2
- e_1}{2}\Bigr)^2 \\
=: & G^2 H^2 - g^2 H^2 - h^2 G^2 + g^2 h^2
\end{align*}
where $G = \frac{d_2 + d_1}{2}$, $H = \frac{e_2 + e_1}{2}$, $g =
\frac{d_2 - d_1}{2}$ and $h = \frac{e_2 - e_1}{2}$. Now we want
\begin{align} \label{gh}
x \approx d_1 e_1 d_2 e_2 = & G^2 H^2 - g^2 H^2 - h^2 G^2 + g^2 h^2 \nonumber \\
G^2 H^2 - x \approx & g^2 H^2 + h^2 G^2 - g^2 h^2 \nonumber \\
(GH - \sqrt{x}) (GH + \sqrt{x}) \approx & g^2 H^2 + h^2 G^2 - g^2
h^2
\end{align}
By $\theta = 1/4$ case in Theorem \ref{thm1}, for some constant $C >
0$, there exist integers $G, H \in [x^{1/4} - C x^{1/16}, x^{1/4} +
C x^{1/16}]$ such that $0 < GH - \sqrt{x} \asymp x^{1/8}$. Then the
left hand side of (\ref{gh}) is $\asymp x^{1/2 + 1/8}$. As for the
right hand side of (\ref{gh}), observe that, for fixed $h$ (say $h =
1$), the increment
$$[(i+1)^2 H^2 + h^2 G^2 - (i+1)^2 h^2] - [i^2 H^2 + h^2 G^2 - i^2
h^2] = (2i + 1) H^2 - (2i + 1) h^2 \asymp x^{1/2} i.$$ Now observe
that
\begin{align*}
& g^2 H^2 + h^2 G^2 - g^2 h^2 \\
=& h^2 G^2 + \sum_{0 \leq i < g} [(i+1)^2 H^2 +
h^2 G^2 - (i+1)^2 h^2] - [i^2 H^2 + h^2 G^2 - i^2 h^2] \\
\asymp & x^{1/2} \sum_{1 \leq i < g} i \asymp g^2 x^{1/2}.
\end{align*}
Therefore, for some integer $1 \leq g \asymp x^{1/16}$,
$$|\hbox{Right hand side of }(\ref{gh}) - \hbox{Left hand side of
}(\ref{gh})| \ll x^{1/2} g \ll x^{1/2 + 1/16}.$$ This gives
$$|x - (G^2 - g^2)(H^2 - h^2)| \ll x^{1/2 + 1/16}$$
or
$$|x - d_1 d_2 e_1 e_2| = |x - (G - g)(G + g)(H - h)(H + h)| \ll
x^{1/2 + 1/16}.$$ Consequently, with
\begin{align*}
a_1 = d_1 e_1 =& (G - g)(H - h) = x^{1/2} + O(x^{1/4 + 1/16}), \\
b_1 = d_2 e_2 =& (G + g)(H + h) = x^{1/2} + O(x^{1/4 + 1/16}), \\
a_2 = d_1 e_2 =& (G - g)(H + h) = x^{1/2} + O(x^{1/4 + 1/16}), \\
b_2 = d_2 e_1 =& (G + g)(H - h) = x^{1/2} + O(x^{1/4 + 1/16}),
\end{align*}
we have a $(\theta,C')$-almost square $n = a_1 b_1 = a_2 b_2$ of
type $2$ in the interval $[x - C'' x^{1/2 + 1/16}, x + C'' x^{1/2 +
1/16}]$ for some $C', C'' > 0$. This proves that $g(\theta) \leq
9/16$ for $\theta \geq 1/4 + 1/16 = 5/16$.
\section{Proof of Theorem \ref{thm3} (iii)}

This time we try to approximate the left hand side of (\ref{gh}) by
the quadratic form $g^2 H^2 + h^2 G^2$ directly. As in the proof of
Theorem \ref{thm3} (iii), for some $C > 0$, there exist integers
$x^{1/4} - C x^{1/16} \leq G, H \leq x^{1/4} + C x^{1/16}$ such that
$0 < GH - \sqrt{x} \asymp x^{1/8}$. The left hand side of (\ref{gh})
is $\asymp x^{1/2 + 1/8}$. Without loss of generality, $G \leq H$.
Then $g^2 H^2 + h^2 G^2 = G^2 (g^2 + h^2) + (H^2 - G^2) g^2$.
Observe that $0 \leq H^2 - G^2 = (H - G)(H + G) \ll x^{1/4 + 1/16}$.
By elementary argument, for any real number $X
> 0$, we can find a sum of two squares $g^2 + h^2$ such that $|X -
(g^2 + h^2)| \ll X^{1/4}$. In particular, we can find $1 \leq g, h
\ll x^{1/16}$ such that
$$\Big|\frac{(GH - \sqrt{x}) (GH + \sqrt{x})}{G^2} - (g^2 + h^2)\Big|
\ll x^{1/32}.$$ This implies
\begin{align*}
& |(GH - \sqrt{x}) (GH + \sqrt{x}) - (g^2 H^2 + h^2 G^2 - g^2 h^2)| \\
\leq & |(GH - \sqrt{x}) (GH + \sqrt{x}) - G^2 (g^2 + h^2)| + |(H^2 -
G^2) g^2| + |g^2 h^2| \ll x^{1/2 + 1/32}.
\end{align*}
Hence
$$|x - d_1 d_2 e_1 e_2| = |x - (G - g)(G + g)(H - h)(H + h)| \ll
x^{1/2 + 1/32}.$$ Consequently, with
\begin{align*}
a_1 = d_1 e_1 =& (G - g)(H - h) = x^{1/2} + O(x^{1/4 + 1/16}), \\
b_1 = d_2 e_2 =& (G + g)(H + h) = x^{1/2} + O(x^{1/4 + 1/16}), \\
a_2 = d_1 e_2 =& (G - g)(H + h) = x^{1/2} + O(x^{1/4 + 1/16}), \\
b_2 = d_2 e_1 =& (G + g)(H - h) = x^{1/2} + O(x^{1/4 + 1/16}),
\end{align*}
there is a $(\theta,C')$-almost square $n = a_1 b_1 = a_2 b_2$ of
type $2$ in the interval $[x - C'' x^{1/2 + 1/32}, x + C'' x^{1/2 +
1/32}]$ for some $C', C'' > 0$. This proves that $g(\theta) \leq
17/32$ for $\theta \geq 1/4 + 1/16 = 5/16$.
\section{Proof of Theorem \ref{thm3} (iv)}

Let $1/2 \leq \phi \leq 1$. Observe that, for large $x$, the
interval $[x + x^{1 - \phi}, x + 2x^{1 - \phi}]$ contains an integer
$n$ which is divisible by an integer $a \in [x^{1 - \phi} / 2, x^{1
- \phi}]$. In particular $n = a b$ with integer $b \in [x^\phi,
3x^\phi]$.

\bigskip

Again we use (\ref{gh}). Instead of having $G, H$ close to $x^{1/4}$
in the proof of Theorem \ref{thm3} (iii), we want
$$G \approx x^{(1 - \phi)/2} \hbox{ and } H \approx x^{\phi/2}
\hbox{ for some } 1/2 < \phi < 2/3.$$ By the observation at the
beginning of this section, we can find $H \in [x^{\phi/2},
3x^{\phi/2}]$ and $G \in [x^{(1-\phi)/2} / 2, x^{(1-\phi)/2}]$ such
that $0 < GH - \sqrt{x} \asymp x^{(1-\phi)/2}$. Then the left hand
side of (\ref{gh}), $L = (GH - \sqrt{x})(GH + \sqrt{x}) \asymp x^{1
- \phi/2}$.

\bigskip

Firstly we approximate $L$ by $g^2 H^2$. For some choice of $g
\asymp x^{1/2 - 3\phi/4}$, we have $0 < L - g^2 H^2 \asymp g H^2
\asymp x^{1/2 + \phi/4}$. Note that $1/2 - 3\phi/4 > 0$ as $\phi <
2/3$.

\bigskip

Secondly we approximate $L - g^2 H^2$ by $h^2 G^2$. For some choice
of $h \asymp x^{5\phi/8 - 1/4}$, we have $|L - g^2 H^2 - h^2 G^2|
\ll h G^2 \asymp x^{3/4 - 3 \phi/8}$. Note that $5 \phi/8 - 1/4
> 0$ as $\phi > 1/2$.

\bigskip

Thirdly, observe that $g^2 h^2 \ll x^{1/2 - \phi/4} \ll x^{3/4 - 3
\phi/8}$ as $\phi < 2$. Therefore, $|L - g^2 H^2 - h^2 G^2 + g^2
h^2| \ll x^{3/4 - 3 \phi/8}$ which gives
$$|x - d_1 d_2 e_1 e_2| = |x - (G - g)(G + g)(H - h)(H + h)| \ll
x^{3/4 - 3 \phi/8}.$$ Consequently, as $1/2 < \phi < 2/3$, with
\begin{align*}
a_1 = d_1 e_1 =& (G - g)(H - h) = x^{1/2} + O(x^{1/2 - \phi/4}), \\
b_1 = d_2 e_2 =& (G + g)(H + h) = x^{1/2} + O(x^{1/2 - \phi/4}), \\
a_2 = d_1 e_2 =& (G - g)(H + h) = x^{1/2} + O(x^{1/2 - \phi/4}), \\
b_2 = d_2 e_1 =& (G + g)(H - h) = x^{1/2} + O(x^{1/2 - \phi/4}),
\end{align*}
there is a $(1/2 - \phi/4, C')$-almost square $n = a_1 b_1 = a_2
b_2$ of type $2$ in the interval $[x - C'' x^{3/4 - 3 \phi/8}, x +
C'' x^{3/4 - 3 \phi/8}]$ for some $C', C'' > 0$. By picking $\phi$
close to $2/3$, we have $g(\theta) \leq 1/2$ for $\theta > 1/3$.
\section{Integer almost divisible by some integer in an interval}

Again let $1/2 \leq \phi \leq 1$. In the previous section, we found
an interval of length $x^{1 - \phi}$ around $x$ containing an
integer divisible by some integer in the interval $[x^{1 - \phi} /
2, x^{1 - \phi}]$. This is obviously true. Our goal in this section
is to find a shorter interval still containing an integer divisible
by some integer in the interval $[x^{1 - \phi} / 2, x^{1 - \phi}]$.
We hope that this will give some improvements to Theorem \ref{thm3}
(iv). Let us reformulate the question as follows:
\begin{ques}
Let $0 < \alpha \leq 1/2$ and $X > 0$ be a large integer. Given $0 <
c_1 < c_2 \leq 1$, find $L$, as small as possible, such that the
interval $[X-L, X]$ contains an integer that is divisible by some
integer in the interval $[c_1 X^\alpha , c_2 X^\alpha]$.
\end{ques}
One may interpret the above as finding an integer in the interval
$[c_1 X^\alpha , c_2 X^\alpha]$ that almost divides $X$ (with a
remainder less than or equal to $L$). We suspect the following
\begin{conj} \label{conjL}
For any $\epsilon > 0$, one can take $L = X^\epsilon$ in the above
question as long as $X$ is sufficiently large in terms of
$\epsilon$.
\end{conj}

However, we can only prove
\begin{prop} \label{prop1}
Suppose $(p,q)$ with $0 \leq p \leq \frac{1}{2} \leq q \leq 1$ is an
exponent pair for exponential sums. Then one can take $L =
X^{\frac{\alpha (q - p)}{1 + p} + \frac{p}{1 + p} + \epsilon}$ in
the above question for any $\epsilon > 0$ as long as $X$ is
sufficiently large in terms of $\epsilon$.
\end{prop}
Our method of proof is making use of Erd\H{o}s-Tur\'{a}n inequality
in the following form (see H.L. Montgomery [\ref{M}, Corollary 1.2]
for example):
\begin{lem} \label{et}
Suppose $M$ is a positive integer chosen so that
$$\sum_{l = 1}^{M} \Big| \sum_{j = 1}^{J} e(l x_j) \Big| \leq
\frac{J}{10}.$$ Then every arc $\mathcal{J} = [\alpha, \beta]
\subseteq [0,1]$ of length $\beta - \alpha \geq \frac{4}{M+1}$
contains at least $\frac{1}{2} J (\beta - \alpha)$ of the points
$x_j$, $1 \leq j \leq J$. Here $||x|| = \min_{n \in \mathbb{Z}} |x -
n|$, the distance from $x$ to the nearest integer, and $e(x) = e^{2
\pi i x}$.
\end{lem}

Proof of Proposition \ref{prop1}: Our sequence $\{x_j\}_{j=1}^{J}$
should be $\{ \frac{X}{a} : a \in \mathbb{Z} \hbox{ and } a \in [c_1
X^\alpha, c_2 X^\alpha] \}$. We want to find some $a$ such that the
fractional part of $\frac{X}{a}$ is small. For if $\{\frac{X}{a}\}
\in [0, \frac{4}{M}]$, then $\frac{X}{a} = k + \frac{\theta}{M}$ for
some integer $k$ and $0 \leq \theta \leq 4$. This gives $X = k a +
\frac{\theta}{M} a$ and $X - \frac{\theta}{M} a = k a$. Hence, with
$L = \frac{4 c_2 X^\alpha}{M}$, the interval $[X-L, X]$ contains an
integer that is divisible by some integer in $[c_1 X^\alpha, c_2
X^\alpha]$. Thus, in view of Lemma \ref{et}, it suffices to show
$$S := \sum_{l = K}^{2K} \Big| \sum_{c_1 X^\alpha \leq a \leq c_2
X^\alpha} e( \frac{l X}{a}) \Big| \leq X^{\alpha - \epsilon}$$ for
any $2 K \leq M$ and $\epsilon > 0$ as long as $X$ is sufficiently
large in terms of $\epsilon$. Keep in mind that we want $M$ as large
as possible.

By the theory of exponent pairs on exponential sums (see Chapter 3
section 4 of [\ref{M}] for an overview),
\begin{equation} \label{exp}
\sum_{c_1 X^\alpha \leq a \leq c_2 X^\alpha} e( \frac{l X}{a}) \ll
(l X (X^\alpha)^{-2} )^p (X^\alpha)^q \ll K^p X^{p - 2 \alpha p +
\alpha q}
\end{equation}
if $(p,q)$ with $0 \leq p \leq \frac{1}{2} \leq q \leq 1$ is an
exponent pair. Using (\ref{exp}), we have
$$\sum_{l = K}^{2K} \Big| \sum_{c_1 X^\alpha \leq a \leq c_2
X^\alpha} e( \frac{l X}{a}) \Big| \ll K^{1+p} X^{p - 2 \alpha p +
\alpha q}.$$ Thus $S \leq X^{\alpha - \epsilon}$ provided that, for
$X$ large enough,
$$K^{1+p} X^{p - 2 \alpha p + \alpha q} \leq X^{\alpha -
\epsilon} \; \hbox{ or } \; K \leq X^{\frac{\alpha (1 - q + 2p)}{1 +
p} - \frac{p}{1 + p} - \epsilon}.$$ Therefore, we can pick
\begin{equation*}
M = X^{\frac{\alpha (1 - q + 2p)}{1 + p} - \frac{p}{1 + p} -
\epsilon} \; \hbox{ which gives } \; L = 4 c_2 X^{\frac{\alpha (q -
p)}{1 + p} + \frac{p}{1 + p} + \epsilon}.
\end{equation*}
This proves Proposition \ref{prop1} since $\epsilon$ is arbitrary.
\section{Proof of Theorem \ref{thm3} (v)}

Proof: We follow closely the proof of the third result. Applying
Proposition \ref{prop1} with $\alpha = 1 - \phi$ and $X = x + 3
x^{\frac{(1 - \phi) (q - p)}{1 + p} + \frac{p}{1 + p} + \epsilon}$,
the interval $[x + x^{\frac{(1 - \phi) (q - p)}{1 + p} + \frac{p}{1
+ p} + \epsilon}, x + 3 x^{\frac{(1 - \phi) (q - p)}{1 + p} +
\frac{p}{1 + p} + \epsilon}]$ contains an integer $n = a b$ with
integers $a \in [x^{1 - \phi} / 2, x^{1 - \phi}]$ and $b \in
[x^\phi, 3 x^\phi]$. Thus we can find
$$H \in [x^{\phi/2}, 3x^{\phi/2}] \hbox{ and } G \in [x^{(1-\phi)/2}
/ 2, x^{(1-\phi)/2}]$$ such that
$$0 < GH - \sqrt{x} \asymp x^{\frac{(1 - \phi) (q - p)}{2 (1 + p)} +
\frac{p}{2(1 + p)} + \frac{\epsilon}{2}}.$$ Then the left hand side
of (\ref{gh}), $L = (GH - \sqrt{x})(GH + \sqrt{x}) \asymp x^{\frac{1
+ p + q}{2 (1 + p)} - \frac{q - p}{2(1 + p)} \phi +
\frac{\epsilon}{2}}$.

\bigskip

Firstly we approximate $L$ by $g^2 H^2$. For some choice of $g
\asymp x^{\frac{1 + p + q}{4(1 + p)} - \frac{2 + p + q}{4(1+p)} \phi
+ \frac{\epsilon}{4}}$, we have $0 < L - g^2 H^2 \asymp g H^2 \asymp
x^{\frac{1 + p + q}{4(1 + p)} + \frac{2 + 3p - q}{4(1+p)} \phi +
\frac{\epsilon}{4}}$. Note that we need $\frac{1 + p + q}{4(1 + p)}
- \frac{2 + p + q}{4(1+p)} \phi \geq 0$ which means $\phi \leq
\frac{1 + p + q}{2 + p + q}$.

\bigskip

Secondly we approximate $L - g^2 H^2$ by $h^2 G^2$. For some choice
of $h \asymp x^{\frac{6 + 7p - q}{8(1+p)} \phi - \frac{3 + 3p -
q}{8(1+p)} + \frac{\epsilon}{8}}$, we have $|L - g^2 H^2 - h^2 G^2|
\ll h G^2 \asymp x^{\frac{5 + 5p + q}{8(1+p)} - \frac{2 + p +
q}{8(1+p)} \phi + \frac{\epsilon}{8}}$. Note that $\frac{6 + 7p -
q}{8(1+p)} \phi - \frac{3 + 3p - q}{8(1+p)} \geq 0$ as $\phi \geq
1/2$ and $p, q \geq 0$.

\bigskip

Thirdly, observe that $g^2 h^2 \ll x^{\frac{3q - p - 1}{4(1+p)} -
\frac{3q - 5p - 2}{4(1+p)} \phi + \frac{3\epsilon}{4}} \ll
x^{\frac{5 + 5p + q}{8(1+p)} - \frac{2 + p + q}{8(1+p)} \phi +
\frac{\epsilon}{8}}$ provided $\phi < \frac{7 + 7p - 5q}{6 + 11p -
5q}$ and $\epsilon$ is small enough. One can easily check that
$\frac{7 + 7p - 5q}{6 + 11p - 5q} > \frac{1 + p + q}{2 + p + q}$.
Therefore, $|L - g^2 H^2 - h^2 G^2 + g^2 h^2| \ll_\epsilon
x^{\frac{5 + 5p + q}{8(1+p)} - \frac{2 + p + q}{8(1+p)} \phi +
\frac{\epsilon}{8}}$ which gives
$$|x - d_1 d_2 e_1 e_2| = |x - (G - g)(G + g)(H - h)(H + h)|
\ll_\epsilon x^{\frac{5 + 5p + q}{8(1+p)} - \frac{2 + p + q}{8(1+p)}
\phi + \frac{\epsilon}{8}}$$ provided $\frac{1}{2} \leq \phi \leq
\frac{1 + p + q}{2 + p + q}$. Choose $\phi = \frac{1 + p + q}{2 + p
+ q}$, we have, after some simple algebra,
$$|x - d_1 d_2 e_1 e_2| = |x - (G - g)(G + g)(H - h)(H + h)|
\ll_\epsilon x^{\frac{1}{2} + \frac{\epsilon}{8}}.$$ Now observe
that with $\phi = \frac{1 + p + q}{2 + p + q}$, after some algebra,
$$GH - x^{1/2} \ll x^{\frac{q}{2(1+p)} - \frac{q - p}{2(1+p)} \phi +
\frac{\epsilon}{2}} = x^{\frac{p+q}{2(2+p+q)} +
\frac{\epsilon}{2}},$$
$$gH \ll x^{\frac{1+p+q}{4(1+p)} - \frac{2+p+q}{4(1+p)} \phi +
\frac{\phi}{2} + \frac{\epsilon}{4}} = x^{\frac{1+p+q}{2(2+p+q)} +
\frac{\epsilon}{4}},$$
$$hG \ll x^{\frac{6 + 7p - q}{8(1+p)} \phi - \frac{3 + 3p -
q}{8(1+p)} + \frac{1 - \phi}{2} + \frac{\epsilon}{8}} =
x^{\frac{1+p+q}{2 (2+p+q)} + \frac{\epsilon}{8}},$$ and
$$gh \ll x^{\frac{3q - p - 1}{8(1+p)} -
\frac{3q - 5p - 2}{8(1+p)} \phi + \frac{3\epsilon}{8}} =
x^{\frac{p+q}{2(2+p+q)} + \frac{3\epsilon}{8}}.$$ Therefore $a_1 =
d_1 e_1 = (G - g)(H - h)$, $b_1 = d_2 e_2 = (G + g)(H + h)$, $a_2 =
d_1 e_2 = (G - g)(H + h)$ and $b_2 = d_2 e_1 = (G + g)(H - h)$ are
all $= x^{\frac{1}{2}} + O_\epsilon(x^{\frac{1+p+q}{2(2+p+q)} +
\frac{\epsilon}{2}})$. Therefore, there is a
$(\frac{1+p+q}{2(2+p+q)} + \frac{\epsilon}{2}, C_\epsilon)$-almost
square of type $2$ in the interval $[x - x^{\frac{1}{2} +
\frac{\epsilon}{8}}, x + x^{\frac{1}{2} + \frac{\epsilon}{8}}]$.
This shows that $g(\theta) \leq 1/2$ for $\theta >
\frac{1+p+q}{2(2+p+q)}$. Since $\frac{1+u}{2+u}$ is an increasing
function, we try to find exponent pairs that make $p+q$ as small as
possible.

For example, recently Huxley [\ref{Hux}] proved that $(p,q) =
(\frac{32}{205} + \epsilon, \frac{1}{2} + \frac{32}{205} +
\epsilon)$ is an exponent pair for any $\epsilon > 0$. This gives
$\frac{1+p+q}{2(2+p+q)} \geq \frac{743}{2306} + \epsilon$ for any
$\epsilon > 0$ and hence Theorem \ref{thm3} (v). Note that
$\frac{743}{2306} = 0.3222029488... < \frac{1}{3}$. However, we
still cannot beat the $\frac{1}{2}$ bound for $g(\theta)$. Assuming
the exponent pair conjecture that $(\epsilon, \frac{1}{2} +
\epsilon)$ is an exponent pair, we can push the range for $\theta$
to $\theta > 0.3$ with $g(\theta) \leq \frac{1}{2}$ but this is
still shy of the range $\theta \geq \frac{1}{4}$. Nevertheless, if
one assumes Conjecture \ref{conjL} in the previous section and
imitates the proof of the (iv) or (v) of Theorem \ref{thm3}, one can
get $g(\theta) \leq \frac{1}{2}$ for $\theta > \frac{1}{4}$. This
comes close to the conjecture $g(1/4) = 1/2$.

\bigskip

{\bf Acknowledgements} The author would like to thank the American
Institute of Mathematics where the study of almost squares began
during a visit from 2004 to 2005. He also thanks Central Michigan
University where the main idea of this paper was worked out during a
one-year visiting position (2005-2006). Finally, he thanks the
University of Hong Kong where the Erd\H{o}s-Tur\'{a}n and exponent
pair part was worked out during a visit there in the summer of 2007.

Department of Mathematical Sciences \\
University of Memphis \\
Memphis, TN 38152 \\
U.S.A. \\
tchan@memphis.edu
\end{document}